\documentclass[a4paper,12pt]{article} 

\usepackage{color}

\usepackage{amsfonts, amsmath, amsthm, amssymb}
\usepackage[T1]{fontenc}
\usepackage[cp1250]{inputenc}
\usepackage{xcolor}
\usepackage{graphicx}
\usepackage{amssymb}
\usepackage{amsmath}
\usepackage{mathptmx}
\usepackage{helvet}
\usepackage{courier}
\usepackage{txfonts}
\usepackage{tikz} 
\usetikzlibrary{arrows}
\usepackage{type1cm}

\usepackage{verbatim}

\usepackage{graphicx}
\usepackage{epsfig,amscd,amssymb,amsxtra,amsmath,amsthm}
\usepackage{type1cm}
\usepackage[T1]{fontenc}
\usepackage{graphics}
\usepackage[mathscr]{eucal}
\usepackage[all]{xy}
\usepackage{amsmath,amscd}

\newtheorem{theorem}{Theorem}[section]
\newtheorem{proposition}[theorem]{Proposition}
\newtheorem{definition}[theorem]{Definition}
\newtheorem{lemma}[theorem]{Lemma}

\newtheorem{remark}[theorem]{Remark}

\newtheorem{observation}[theorem]{Observation}

\begin{document}

\def\joinrel{\mkern-3mu}
\newcommand{\varproj}{\displaystyle \lim_{\multimapinv\joinrel-\joinrel-}}

\title{An uncountable family of non-smooth fans  that admit transitive homeomorphisms}
\author{Iztok Bani\v c, Judy Kennedy,  Chris Mouron and Van Nall}
\date{}

\maketitle

\begin{abstract}
Recently, many examples of smooth fans that admit a transitive homeomorphism have been constructed. For example, a family of uncountably many pairwise non-homeomorphic smooth fans that admit transitive homeomorphisms was constructed in \cite{banic1}. In this paper, we construct a family of uncountably many pairwise non-homeomorphic non-smooth fans that admit transitive homeomorphisms. 
\end{abstract}
\-

\section{Introduction}

Recently, many examples of smooth fans that admit a transitive homeomorphism have been constructed. For example, in \cite{BE}, a transitive homeomorphism was constructed on the Cantor fan, in \cite{banic1}, a transitive homeomorphism was constructed on the star of Cantor fans, and in \cite{banic2} and \cite{oprocha}, a transitive homeomorphism was constructed on the Lelek fan. The homeomorphism from \cite{oprocha} has 0-entropy while the homeomorphism from \cite{banic2} has positive entropy. Then, an uncountable family of smooth fans that admit  transitive homeomorphisms was constructed in \cite{banic1}. In this paper, we show that there are also examples of non-smooth fans that admit transitive homeomorphisms. Even more, we show that there is a family of uncountably many pairwise non-homeomorphic non-smooth fans that admit transitive homeomorphisms.

We proceed as follows. In Section \ref{s1}, we introduce the definitions, notation and the well-known results that will be used later in the paper and in Section \ref{s2}, a family of uncountably many pairwise non-homeomorphic non-smooth fans that admit a transitive homeomorphism is constructed.

\section{Definitions and Notation}\label{s1}
The following definitions, notation and well-known results are needed in the paper.


\begin{definition}
Let $(X,d)$ be a metric space, $x\in X$ and $\varepsilon>0$. We use $B(x,\varepsilon)$ to denote the open ball,  {centered} at $x$ with radius $\varepsilon$.
\end{definition}
\begin{definition}
Let $(X,d)$ be a compact metric space. Then we define \emph{$2^X$} by 
$$
2^{X}=\{A\subseteq X \ | \ A \textup{ is a non-empty closed subset of } X\}.
$$
Let $\varepsilon >0$ and let $A\in 2^X$. Then we define  \emph{$N_d(\varepsilon,A)$} by 
$$
N_d(\varepsilon,A)=\bigcup_{a\in A}B(a,\varepsilon).
$$
Let $A,B\in 2^X$. The function \emph{$H_d:2^X\times 2^X\rightarrow \mathbb R$}, defined by
$$
H_d(A,B)=\inf\{\varepsilon>0 \ | \ A\subseteq N_d(\varepsilon,B), B\subseteq N_d(\varepsilon,A)\},
$$
is called \emph{the Hausdorff metric}. The Hausdorff metric is in fact a metric and the metric space $(2^X,H_d)$ is called \emph{a hyperspace of the space $(X,d)$}. 
\end{definition}
\begin{remark}
Let $(X,d)$ be a compact metric space, let $A$ be a non-empty closed subset of $X$,  and let $(A_n)$ be a sequence of non-empty closed subsets of $X$. When we say $\displaystyle A=\lim_{n\to \infty}A_n$ with respect to the {Hausdorff} metric, we mean $\displaystyle A=\lim_{n\to \infty}A_n$ in $(2^X,H_d)$. 
\end{remark}
\begin{definition}
 \emph{A continuum} is a non-empty compact connected metric space.  \emph{A subcontinuum} is a subspace of a continuum, which is itself a continuum.
 \end{definition}
 \begin{definition}
Let $X$ be a continuum. 
\begin{enumerate}
\item The continuum $X$ is \emph{ unicoherent}, if for any subcontinua $A$ and $B$ of $X$ such that $X=A\cup B$,  the compactum $A\cap B$ is connected. 
\item The continuum $X$ is \emph{hereditarily unicoherent } provided that each of its subcontinua is unicoherent.
\item The continuum $X$ is a \emph{dendroid}, if it is an arcwise connected, hereditarily unicoherent continuum.
\item Let $X$ be a continuum.  If $X$ is homeomorphic to $[0,1]$, then $X$ is \emph{ an arc}.   
\item A point $x$ in an arc $X$ is called \emph{an end-point of the arc  $X$}, if  there is a homeomorphism $\varphi:[0,1]\rightarrow X$ such that $\varphi(0)=x$.
\item Let $X$ be a dendroid.  A point $x\in X$ is called an \emph{end-point of the dendroid $X$}, if for  every arc $A$ in $X$ that contains $x$, $x$ is an end-point of $A$.  The set of all end-points of $X$ will be denoted by $E(X)$. 
\item A continuum $X$ is \emph{a simple triod}, if it is homeomorphic to $([-1,1]\times \{0\})\cup (\{0\}\times [0,1])$.
\item A point $x$ in a simple triod $X$ is called \emph{the top-point} or just the \emph{top of the simple triod $X$}, if  there is a homeomorphism $\varphi:([-1,1]\times \{0\})\cup (\{0\}\times [0,1])\rightarrow X$ such that $\varphi(0,0)=x$.
\item Let $X$ be a dendroid.  A point $x\in X$ is called \emph{a ramification-point of the dendroid $X$}, if there is a simple triod $T$ contained in $X$ such that $x$ is the top of $T$.  The set of all ramification-points of $X$ will be denoted by $R(X)$. 
\item The continuum $X$ is \emph{a  fan}, if it is a dendroid with at most one ramification point $v$, which is called the top of the fan $X$ (if it exists).
\item Let $X$ be a fan.   For all points $x$ and $y$ in $X$, we define  \emph{$A_X[x,y]$} to be the arc in $X$ with end-points $x$ and $y$, if $x\neq y$. If $x=y$, then we define $A_X[x,y]=\{x\}$.
\item Let $X$ be a fan with top $v$. We say that that the fan $X$ is \emph{smooth} if for any $x\in X$ and for any sequence $(x_n)$ of points in $X$,
$$
\lim_{n\to \infty}x_n=x \Longrightarrow \lim_{n\to \infty}A_X[v,x_n]=A_X[v,x]
$$ 
with respect to the Hausdorff metric.
\item Let $X$ be a fan.  We say that $X$ is \emph{a Cantor fan}, if $X$ is homeomorphic to the continuum
$$
\bigcup_{c\in C}{S}_c,
$$
where $C\subseteq [0,1]$ is the standard Cantor set and for each $c\in C$, ${S}_c$ is the straight line segment in the plane from $(0,0)$ to $(c,1)$. See Figure \ref{fig000}, where a Cantor fan is pictured.
\begin{figure}[h!]
	\centering
		\includegraphics[width=30em]{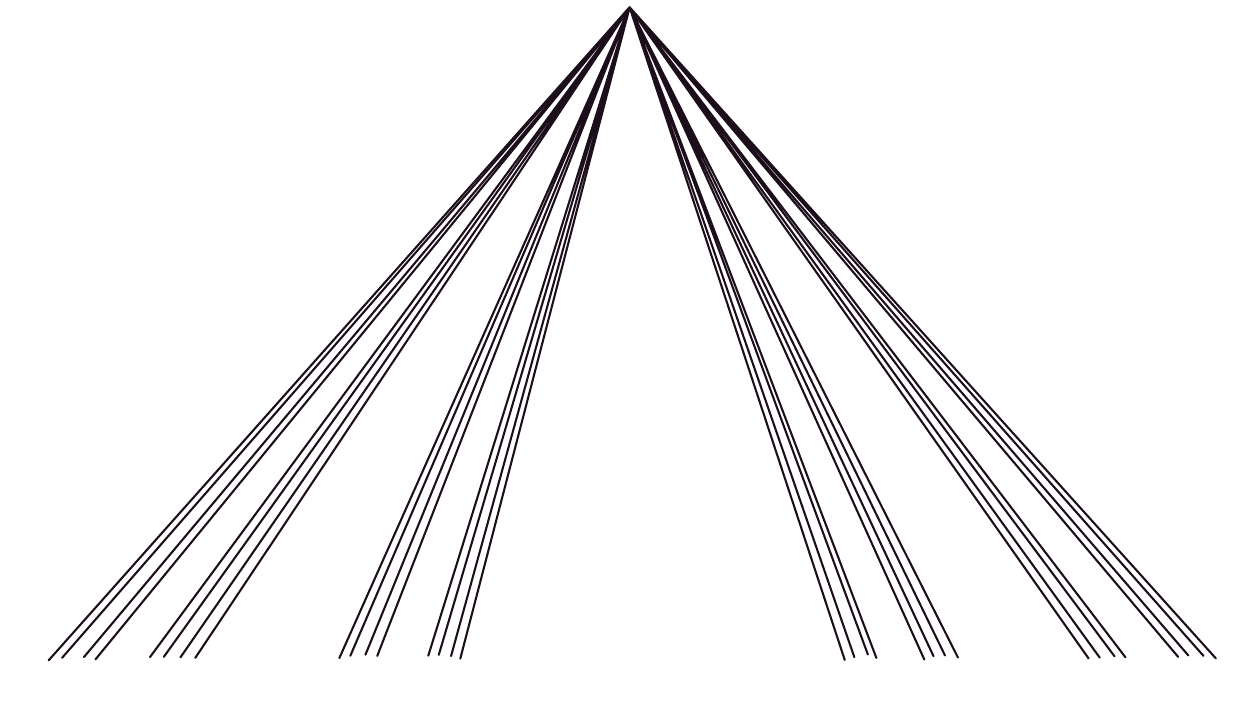}
	\caption{The Cantor fan}
	\label{fig000}
\end{figure}  
\end{enumerate}
\end{definition}
\begin{observation}\label{embi}
	It is a well-known fact that the Cantor fan is universal for smooth fans, {i.e., every smooth fan embeds into the Cantor fan} (for details see \cite[Theorem 9, p. 27]{Jcharatonik},  \cite[ Corollary 4]{koch},  and \cite{eberhart}).
\end{observation}
\begin{definition}
	Let $F$ be a fan with top $o$. For each end-point $e\in E(F)$ of the fan $F$, the arc in $F$ from $o$ to $e$ is called \emph{a leg of $F$}. The set of all legs of $F$ is denoted by $\mathcal L(F)$.
\end{definition}
\begin{definition}
Let $(X,f)$ be a dynamical system.  We say that $(X,f)$ is  \emph{transitive}, if for all non-empty open sets $U$ and $V$ in $X$,  there is a non-negative integer $n$ such that $f^n(U)\cap V\neq \emptyset$. We say that the mapping $f$ is \emph{transitive}, if $(X,f)$ is transitive.
\end{definition}
\begin{definition}\label{povezava}
Let $X$ be a compact metric space. We say that $X$  \emph{admits a transitive homeomorphism}, if there is a homeomorphism $f:X\rightarrow X$ such that $(X,f)$ is transitive.
\end{definition}
\begin{definition}
Let $X$ be a compact metric space and let $\sim$ be an equivalence relation on $X$. For each $x\in X$, we use $[x]$ to denote the equivalence class of the element $x$ with respect to the relation $\sim$. We also use $X/_{\sim}$ to denote the quotient space $X/_{\sim}=\{[x] \ | \ x\in X\}$. 
\end{definition}
\begin{observation}\label{kvokvo}
	Let $X$ be a compact metric space, let $\sim$ be an equivalence relation on $X$, let $q:X\rightarrow X/_{\sim}$ be the quotient map that is defined by $q(x)=[x]$ for each $x\in X$,  and let $U\subseteq X/_{\sim}$. Then 
	$$
	U \textup{ is open in } X/_{\sim} ~~~ \Longleftrightarrow ~~~ q^{-1}(U)\textup{ is open in } X.
	$$
\end{observation}
\begin{definition}
Let $X$ be a compact metric space, let $\sim$ be an equivalence relation on $X$,  and let $f:X\rightarrow X$ be a function such that for all $x,y\in X$,
$$
x\sim y  \Longleftrightarrow f(x)\sim f(y).
$$
 Then we let $f^{\star}:X/_{\sim}\rightarrow X/_{\sim}$ be defined by   
$
f^{\star}([x])=[f(x)]
$
for any $x\in X$. 
\end{definition}
We also use the following proposition, which is a well-known result.
\begin{proposition}\label{kvocienti}
Let $X$ be a compact metric space, let $\sim$ be an equivalence relation on $X$, and  let $f:X\rightarrow X$ be a function such that for all $x,y\in X$,
$$
x\sim y  \Longleftrightarrow f(x)\sim f(y).
$$
 Then the following hold.
\begin{enumerate}
\item $f^{\star}$ is a well-defined function from  $X/_{\sim}$ to $X/_{\sim}$. 
\item If $f$ is continuous, then $f^{\star}$ is continuous.
\item If $f$ is a homeomorphism, then $f^{\star}$ is a homeomorphism.
\item If $f$ is transitive, then $f^{\star}$ is transitive.
\end{enumerate}
\end{proposition}
\begin{proof}
	See \cite[Theorem 3.4]{BE}. 
\end{proof}
\section{Non-smooth fans that admit transitive homeomorphisms}\label{s2}
In this section, we construct a family of uncountably many pairwise non-homeo\-morphic non-smooth fans that admit transitive homeomorphisms. In \cite{banic1}, the notion of $\textup{JuMa}(X)$ for a fan $X$ is introduced. We also use this notion here, i.e., we use Theorems \ref{juma} and \ref{jumas} that are proved in \cite{banic1}. 
\begin{definition}
	Let $X$ be a fan with the top $o$. We define the set $\textup{JuMa}(X)$ as follows: 
	$$
	\textup{JuMa}(X)=\{x\in X\setminus \{o\} \ | \ \textup{ there is a sequence } (e_n) \textup{ in } E(X) \textup{ such that } \lim_{n\to \infty}e_n=x\}.  
$$
\end{definition} 

\begin{theorem}\label{juma}\cite[Proposition 4.53]{banic1}
	Let $X$ and $Y$ be fans with tops $o_X$ and $o_Y$, respectively,  and let $f:X\rightarrow Y$ be a homeomorphism. Then for each $e\in E(X)$, 
	$$
	|A_X[o_X,e]\cap \textup{JuMa}(X)|=|A_Y[o_Y,f(e)]\cap\textup{JuMa}(Y)|.
	$$
	Here $|S|$ denotes the cardinality of $S$ for any set $S$.
\end{theorem}

\begin{theorem}\label{jumas} \cite[Corollary 4.54]{banic1}
	Let $X$ and $Y$ be fans with tops $o_X$ and $o_Y$, respectively. If there is $e\in E(X)$ such that for each $f\in E(Y)$,
	$$
	 |A_Y[o_Y,f]\cap\textup{JuMa}(Y)|\neq |A_X[o_X,e]\cap \textup{JuMa}(X)|,
	$$
	then $X$ and $Y$ are not homeomorphic. 
\end{theorem}

In \cite[Definition 4.50]{banic1}, a family $\mathcal F$ of uncountably many pairwise non-homeo\-morphic smooth fans that admit transitive homeomorphisms was constructed. The family $\mathcal F$ is constructed in such a way that for any two smooth fans $X,Y\in \mathcal F$ with tops $o_X$ and $o_Y$, respectively, there is $e\in E(X)$ such that for each $f\in E(Y)$,
	$$
	 |A_Y[o_Y,f]\cap\textup{JuMa}(Y)|\neq |A_X[o_X,e]\cap \textup{JuMa}(X)|.
	$$
    Theorem \ref{jumas} then implies that $X$ and $Y$ are not homeomorphic. It also follows from the construction of the family $\mathcal F$ that for each $X\in \mathcal F$ there is a transitive homeomorphism $\varphi:X\rightarrow X$ and a leg $A\in \mathcal L(X)$ such that 
    \begin{enumerate}
       	\item for each $x\in A$, $\varphi(x)=x$, and
    	\item $A\cap\textup{JuMa}(X)=\emptyset$.
    \end{enumerate}
 Since all the fans in the family $\mathcal F$ from \cite[Definition 4.50]{banic1} are smooth, they can all be embedded into the Cantor fan as seen by Observation \ref{embi}.
\begin{definition}
	Let $F$ be the Cantor fan, defined by $F=\bigcup_{c\in C}{S}_c$, 
where $C\subseteq [0,1]$ is the standard Cantor set and for each $c\in C$, ${S}_c$ is the straight line segment in the plane from $(0,0)$ to $(c,1)$. For each $X \in \mathcal F$, let $\Psi_{X}:X\rightarrow F$ be an embedding. We use $\mathcal E$ to denote the family
	$$
	\mathcal E=\{\Psi_X(X) \ |\ X\in \mathcal F\}.
	$$
	We  denote the members of the family $\mathcal E$ by $F_{\lambda}$:
	$$
	\mathcal E=\{F_{\lambda} \ | \ \lambda \in \Lambda\}.
	$$
	\end{definition}
\begin{observation}
	Note that for each $\lambda \in \Lambda$,  
	 $F_{\lambda}$ is a smooth fan such that $F_{\lambda}\subseteq F$. Also, for each $\lambda \in \Lambda$, 
	 	\begin{enumerate}
		\item   there is a transitive homeomorphisms $\varphi_{\lambda}:F_{\lambda}\rightarrow F_{\lambda}$ 
		\item there is a leg $A_{\lambda}\in \mathcal L(F_{\lambda})$ 
	\end{enumerate} 
	such that 
		\begin{enumerate}
		\item for each $x\in A_{\lambda}$, $\varphi_{\lambda}(x)=x$, and
    	\item $A_{\lambda}\cap\textup{JuMa}(F_{\lambda})=\emptyset$.
    \end{enumerate}
 For each $\lambda \in \Lambda$, we choose and fix such a homeomorphism $\varphi_{\lambda}$ and such a leg $A_{\lambda}$. We also assume for the rest of the paper that for all $\lambda_1,\lambda_2\in \Lambda$,
 $$
 \lambda_1\neq \lambda_2 ~~~ \Longrightarrow ~~~ F_{\lambda_1} \textup{ is not homeomorphic to } F_{\lambda_2}. 
 $$
\end{observation}
\begin{definition}
	For each $\lambda \in \Lambda$, we use $a_\lambda$ to denote the end-point of $F_{\lambda}$ that is defined as follows:
	$$
E(F_{\lambda})\cap A_{\lambda}=\{a_{\lambda}\}.
	$$
\end{definition}
\begin{observation}
	Note that for each $\lambda\in \Lambda$,
	$$
	A_{\lambda}=\{t\cdot a_{\lambda} \ | \ t\in [0,1]\}.
	$$
\end{observation}
\begin{definition}
	For each $\lambda \in \Lambda$,  we define the relation $\sim_{\lambda}$ on $F_{\lambda}$ as follows. For each $\lambda \in \Lambda$ and for all $x,y\in F_{\lambda}$, we define that 
	$$
	x\sim_{\lambda}y ~~~ \Longleftrightarrow ~~~ x=y \textup{ or there is } t\in [0,1] \textup{ such that } x=t\cdot a_{\lambda} \textup{ and } y=(1-t)\cdot a_{\lambda}.
	$$ 
\end{definition}
\begin{proposition}\label{itak}
	For each $\lambda\in \Lambda$, $\sim_{\lambda}$ is an equivalence relation on $F_{\lambda}$ such that for all $x,y\in F_{\lambda}$,
$$
x\sim_{\lambda} y  \Longleftrightarrow \varphi_{\lambda}(x)\sim_{\lambda} \varphi_{\lambda}(y).
$$
\end{proposition}
\begin{proof}
Let $\lambda\in \Lambda$. It follows from the definition of $\sim_{\lambda}$ that it is an equivalence relation on $F_{\lambda}$. Next, let $x,y\in F_{\lambda}$ be such that $x\sim_{\lambda} y$. We consider the following cases.
\begin{enumerate}
	\item $x=y$. Then $\varphi_{\lambda}(x)=\varphi_{\lambda}(y)$ and it follows that $\varphi_{\lambda}(x)\sim_{\lambda} \varphi_{\lambda}(y)$. 
	\item $x\neq y$. Then $x,y\in A_{\lambda}$. Therefore, $\varphi_{\lambda}(x)=x$ and $\varphi_{\lambda}(y)=y$. It follows that $\varphi_{\lambda}(x)\sim_{\lambda} \varphi_{\lambda}(y)$. 
\end{enumerate}
Finally, let $x,y\in F_{\lambda}$ be such that $\varphi_{\lambda}(x)\sim_{\lambda} \varphi_{\lambda}(y)$. We consider the following cases.
\begin{enumerate}
	\item $\varphi_{\lambda}(x)=\varphi_{\lambda}(y)$. Then $x=y$ (since $\varphi_{\lambda}$ is a homeomorphism) and it follows that $x\sim_{\lambda} y$. 
	\item $\varphi_{\lambda}(x)\neq\varphi_{\lambda}(y)$. Then $\varphi_{\lambda}(x),\varphi_{\lambda}(y)\in A_{\lambda}$. Therefore, $\varphi_{\lambda}(x)=x$ and $\varphi_{\lambda}(y)=y$. It follows that $x\sim_{\lambda} y$. 
\end{enumerate}
\end{proof}	
\begin{definition}
	For each $\lambda\in \Lambda$ and for each $x\in F_{\lambda}$, we define $[x]_{\lambda}$ to be the equivalence class of $x$ with respect to the relation $\sim_{\lambda}$:
	$$
	[x]_{\lambda}=\{y\in F_{\lambda} \ | \ y\sim_{\lambda} x\}.
	$$
	We also define $J_{\lambda}$ to be the quotient space 
	$$
	J_{\lambda}=F_{\lambda}/_{\sim_{\lambda}}
	$$
	and we use $q_{\lambda}$ to denote the quotient map $q_{\lambda}:F_{\lambda}\rightarrow J_{\lambda}$, defined by 
	 $$
 q_{\lambda}(x)=[x]_{\lambda}
 $$
 for each $x\in F_{\lambda}$. 
\end{definition}
\begin{observation}
	It follows from \cite[Theorem 4.2.13]{engelking1} that for each $\lambda\in \Lambda$, $J_{\lambda}$ is metrizable. Since for each $\lambda\in \Lambda$, $F_{\lambda}$ is connected and compact and since $q_{\lambda}$ is continuous, it follows from $J_{\lambda}=q_{\lambda}(F_{\lambda})$ that also $J_{\lambda}$ is connected and compact. Therefore,  for each $\lambda\in \Lambda$, $J_{\lambda}$ is a continuum. 
\end{observation}
\begin{theorem}\label{trf}
	For each $\lambda\in \Lambda$, $\varphi_{\lambda}^{\star}:J_{\lambda}\rightarrow J_{\lambda}$, defined by $\varphi_{\lambda}^{\star}([x]_{\lambda})=[\varphi_{\lambda}(x)]_{\lambda}$ for each $x\in F_{\lambda}$, is a transitive homeomorphism. 
\end{theorem}
\begin{proof}
	For each $\lambda\in \Lambda$, $\varphi_{\lambda}$ is (by Proposition \ref{itak}) a transitive homeomorphism, such that for all $x,y\in F_{\lambda}$,
$$
x\sim_{\lambda} y  \Longleftrightarrow \varphi_{\lambda}(x)\sim_{\lambda} \varphi_{\lambda}(y).
$$
It follows from Proposition \ref{kvocienti} that for each $\lambda\in \Lambda$, $\varphi_{\lambda}^{\star}$ is a transitive homeomorphism.
\end{proof}
\begin{definition}
	For each $\lambda \in \Lambda$ and for each $e\in E(F_{\lambda})$, we define the subsets $K_e\subseteq F_{\lambda}$ and $L_e\subseteq J_{\lambda}$ by
$$
K_e=\{e\cdot t \ | \ t\in [0,1]\} ~~~ \textup{ and } ~~~ L_e=\{[e\cdot t]_{\lambda} \ | \ t\in [0,1]\}. 
$$
\end{definition}
\begin{observation}
	Let $\lambda \in \Lambda$. Note that for each $e\in E(F_{\lambda})$, $K_e\in \mathcal L(F_{\lambda})$, and that  $K_{a_{\lambda}}=A_{\lambda}$.
\end{observation}
\begin{proposition}\label{prop}
For each $\lambda \in \Lambda$, and for each  $e\in E(F_{\lambda})$, $L_e$ is an arc in $J_{\lambda}$. 
\end{proposition}
\begin{proof}
Let $\lambda \in \Lambda$ and let  $e\in E(F_{\lambda})$. We consider the following cases.
	\begin{enumerate}
	\item $e=a_{\lambda}$. Note that
	$$ 
	L_{a_{\lambda}}=\Big\{[a_{\lambda}\cdot t]_{\lambda} \ | \ t\in \Big[0,\frac{1}{2}\Big]\Big\}=\Big\{\{t\cdot a_{\lambda},(1-t)\cdot a_{\lambda}\} \ | \ t\in \Big[0,\frac{1}{2}\Big]\Big\}.
	$$
	Let $h:L_{a_{\lambda}}\rightarrow [0,1]$ be defined by 
	$$
	h([a_{\lambda}\cdot t]_{\lambda})=2t
	$$ 
	for each $t\in [0,\frac{1}{2}]$. We show that $h$ is a homeomorphism. Note that $h$ is a bijection. To show that $h$ is continuous, let $U$ be any open set in $[0,1]$. We show that $h^{-1}(U)$ is open in $L_{a_{\lambda}}$ by using Observation \ref{kvokvo}. Note that 
        $$
        h^{-1}(U)=\Big\{\Big[\frac{1}{2}t\cdot a_{\lambda}\Big]_{\lambda} \ | \ t\in U\Big\}=\Big\{\Big\{\frac{1}{2}t\cdot a_{\lambda},\Big(1-\frac{1}{2}t\Big)\cdot a_{\lambda}\Big\} \ | \ t\in U\Big\}
        $$
        and that
        $$
        q_{\lambda}^{-1}(h^{-1}(U))=\Big\{\frac{1}{2}t\cdot a_{\lambda} \ | \ t\in U\Big\}\cup \Big\{\Big(1-\frac{1}{2}t\Big)\cdot a_{\lambda} \ | \ t\in U\Big\}.
        $$
        Then $q_{\lambda}^{-1}(h^{-1}(U))$ is a union of two open sets in $A_{\lambda}$, therefore, $q_{\lambda}^{-1}(h^{-1}(U))$ is open in $A_{\lambda}$. It follows that $h^{-1}(U)$ is open in $L_{a_{\lambda}}$. This proves that $h$ is continuous. Note that $A_{\lambda}$ is compact (since it is an arc), therefore, since $L_{a_{\lambda}}=q_{\lambda}(A_{\lambda})$, it follows that $L_{a_{\lambda}}$ is compact. Hence, $h$ is a continuous surjection from a compact space to a metric space. It follows that $h$ is a homeomorphism. This proves that $L_{a_{\lambda}}$ is an arc in $J_{\lambda}$. 
    
	\item $e\neq a_{\lambda}$. We show that $L_{e}$ is an arc by showing that it is homeomorphic to $K_e$. Let $H:L_{e}\rightarrow K_{e}$ be defined by 
	$$
	H([t\cdot e]_{\lambda})=t\cdot e
	$$ 
	for each $t\in [0,1]$. We show that $H$ is a homeomorphism. Note that $H$ is a bijection and that for each $t\in [0,1]$, $[t\cdot e]_{\lambda}=\{t\cdot e\}$. Next, we show that $H$ is continuous. Let $U$ be any open set in $K_{e}$. 
	Note that $q_{\lambda}^{-1}(H^{-1}(U))=U$. It follows that $q_{\lambda}^{-1}(H^{-1}(U))$ is open in $K_{e}$. Therefore, $H^{-1}(U)$ is open in $L_{e}$ by Observation \ref{kvokvo}. This proves that $H$ is continuous. Note that $K_{e}$ is compact (since it is an arc), therefore, since $L_{e}=q_{\lambda}(K_{e})$, it follows that $L_{e}$ is compact. Hence, $H$ is a continuous surjection from a compact space to a metric space. It follows that $H$ is a homeomorphism. 	Therefore, $L_e$ is an arc in $J_{\lambda}$. 
\end{enumerate}
Therefore, for each  $e\in E(F_{\lambda})$, $L_e$ is an arc in $J_{\lambda}$. 
\end{proof}
\begin{observation}
	Let $\lambda \in \Lambda$. Note that $L_{a_{\lambda}}$ is an arc with end-points $[(0,0)]_{\lambda}$ and $[\frac{1}{2}\cdot a_{\lambda}]_{\lambda}$.
\end{observation}

In Theorem \ref{fanci}, we prove that each $J_{\lambda}$ is a fan. In its proof, we use Lemma \ref{ema}.
\begin{lemma}\label{ema}
	Let $\lambda\in \Lambda$ and let $A$ be a subcontinuum of $J_{\lambda}$. The following statements are equivalent.
	\begin{enumerate}
		\item\label{m1} The preimage $q_{\lambda}^{-1}(A)$ is not connected.
		\item\label{m2} $A\cap L_{a_\lambda}=\{[(0,0)]_{\lambda}\}$ or there is $s\in \Big(0,\frac{1}{2}\Big)$ such that $A\cap L_{a_\lambda}=\{[t\cdot a_{\lambda}]_{\lambda} \ | \ t\in [0,s]\}$.
	\end{enumerate} 
	\end{lemma}
\begin{proof}
To prove the implication from \ref{m1} to \ref{m2}, suppose that $A\cap L_{a_\lambda}\neq \{[(0,0)]_{\lambda}\}$ and that for each $s\in \Big(0,\frac{1}{2}\Big)$,  $A\cap L_{a_\lambda}\neq \{[t\cdot a_{\lambda}]_{\lambda} \ | \ t\in [0,s]\}$. Then $A\cap L_{a_\lambda}=\emptyset$ or $A\cap L_{a_\lambda}=L_{a_\lambda}$. In both cases, $q_{\lambda}^{-1}(A)$ is connected. Next, we prove the implication from \ref{m2} to \ref{m1}. If $A\cap L_{\lambda}=\{[(0,0)]_{\lambda}\}$, then $a_{\lambda}$ is an isolated point of $q_{\lambda}^{-1}(A)$. Therefore, in this case, $q_{\lambda}^{-1}(A)$ is not connected. Next, suppose that there is $s\in \Big(0,\frac{1}{2}\Big)$ such that $A\cap L_{\lambda}=\{[t\cdot a_{\lambda}]_{\lambda} \ | \ t\in [0,s]\}$.
Choose and fix such an $s$. Let $U=\{(1-t)\cdot a_{\lambda} \ | \ t\in [0,s]\}$. Then $U$ is clopen in $q_{\lambda}^{-1}(A)$. Since $U\neq q_{\lambda}^{-1}(A)$, it follows that also in this case, $q_{\lambda}^{-1}(A)$ is not connected. 
\end{proof}
\begin{observation}\label{preimage}
	Let $\lambda\in \Lambda$ and let $A$ be a subcontinuum of $J_{\lambda}$ such that the preimage $q_{\lambda}^{-1}(A)$ is not connected. Note that
	\begin{enumerate}
		\item if $A\cap L_{a_\lambda}=\{[(0,0)]_{\lambda}\}$, then $q_{\lambda}^{-1}(A)$ has exactly two components, $C_1=\{a_{\lambda}\}$ and $C_2=q_{\lambda}^{-1}(A)\setminus C_1$.
		\item if there is $s\in \Big(0,\frac{1}{2}\Big)$ such that $A\cap L_{a_\lambda}=\{[t\cdot a_{\lambda}]_{\lambda} \ | \ t\in [0,s]\}$, then $q_{\lambda}^{-1}(A)$ has exactly two components, $C_1=\{(1-t)\cdot a_{\lambda} \ | \ t\in [0,s]\}$ and $C_2=q_{\lambda}^{-1}(A)\setminus C_1$.
	\end{enumerate}
\end{observation}
\begin{definition}
	Let $X$ and $Y$ be continua and let $f:X\rightarrow Y$ be a continuous function. We say that $f$ is semi-confluent, if for any subcontinuum $C$ of $Y$ and for all components $C_1$ and $C_2$ of $f^{-1}(C)$, $f(C_1)\subseteq f(C_2)$ or $f(C_2)\subseteq f(C_1)$.   
\end{definition}

\begin{theorem}\label{fanci}
	For each $\lambda\in \Lambda$, $J_{\lambda}$ is a fan.
\end{theorem}
\begin{proof}
Let $\lambda\in \Lambda$. We show that the quotient map $q_{\lambda}:F_{\lambda}\rightarrow J_{\lambda}$ is a semi-confluent surjection. Since $q_{\lambda}$ is a quotient map, it is a surjection. To see that it is semi-confluent, let $C$ be any subcontinuum of $J_{\lambda}$. We consider the following cases.
\begin{enumerate}
	\item $q_{\lambda}^{-1}(C)$ is connected. Then for all components $C_1$ and $C_2$ of $q_{\lambda}^{-1}(C)$, $q_{\lambda}(C_1)\subseteq q_{\lambda}(C_2)$ or $q_{\lambda}(C_2)\subseteq q_{\lambda}(C_1)$ (since $C_1=C_2=q_{\lambda}^{-1}(C)$).   
	\item $q_{\lambda}^{-1}(C)$ is not connected. By Lemma \ref{ema}, $q_{\lambda}^{-1}(C)\cap L_{a_\lambda}=\{[(0,0)]_{\lambda}\}$ or there is $s\in \Big(0,\frac{1}{2}\Big)$ such that $q_{\lambda}^{-1}(C)\cap L_{a_\lambda}=\{[t\cdot a_{\lambda}]_{\lambda} \ | \ t\in [0,s]\}$. According to this, we consider the following two cases.
	\begin{enumerate}
	\item $q_{\lambda}^{-1}(C)\cap L_{a_\lambda}=\{[(0,0)]_{\lambda}\}$. By Observation \ref{preimage}, $q_{\lambda}^{-1}(C)$ has exactly two components: $C_1=\{a_{\lambda}\}$ and $C_2=q_{\lambda}^{-1}(C)\setminus C_1$. Note that $q_{\lambda}(C_1)\subseteq q_{\lambda}(C_2)$.
	\item There is $s\in \Big(0,\frac{1}{2}\Big)$ such that $q_{\lambda}^{-1}(C)\cap L_{a_\lambda}=\{[t\cdot a_{\lambda}]_{\lambda} \ | \ t\in [0,s]\}$. By Observation \ref{preimage}, $q_{\lambda}^{-1}(C)$ has exactly two components, $C_1=\{(1-t)\cdot a_{\lambda} \ | \ t\in [0,s]\}$ and $C_2=q_{\lambda}^{-1}(C)\setminus C_1$. Note that also in this case, $q_{\lambda}(C_1)\subseteq q_{\lambda}(C_2)$.
	\end{enumerate}
\end{enumerate} 
It follows that $q_{\lambda}:F_{\lambda}\rightarrow J_{\lambda}$ is a semi-confluent surjection.  Now it follows from  Ma\' ckowiak's result \cite[Theorem 5.6]{mackowjak}, which says that a semi-confluent image of a fan is a fan, that $J_{\lambda}$ is a fan.
\end{proof}
\begin{theorem}\label{smoothf}
	For each $\lambda\in \Lambda$, the fan $J_{\lambda}$ is not smooth.
\end{theorem}
\begin{proof}
Let $\lambda\in \Lambda$ and let $(x_k)$ be a sequence of points in $F_{\lambda}$ such that
\begin{enumerate}
    \item for each positive integer $k$, there is $e_k\in E(F_{\lambda})\setminus\{a_{\lambda}\}$ such that 
    $$
    x_k\in K_{e_k}\setminus \{(0,0),e_k\},
    $$
    \item for all positive integers $k,\ell$ and for all $e,f\in E(F_{\lambda})$ such that $x_k\in K_{e}\setminus \{(0,0)\}$ and $x_{\ell}\in K_{f}\setminus \{(0,0)\}$,
    $$
    k\neq \ell  ~~~ \Longrightarrow ~~~ e\neq f.  
    $$
	\item $\displaystyle\lim_{k\to \infty}x_k=\frac{3}{4}\cdot a_{\lambda}$.
\end{enumerate} 
It follows from the construction of the family $\mathcal F$ in \cite[Definition 4.50]{banic1} that such a sequence exists. Since $q_{\lambda}$ is continuous, it follows that 
$$
\lim_{k\to \infty}q_{\lambda}(x_k)=q_{\lambda}\Big(\frac{3}{4}\cdot a_{\lambda}\Big)
$$
Note that for each positive integer $k$, 
	$$
	q_{\lambda}(x_k)=[x_k]_{\lambda}=\{x_{\lambda}\}
	$$
and that 
$$
q_{\lambda}\Big(\frac{3}{4}\cdot a_{\lambda}\Big)=\Big[\frac{3}{4}\cdot a_{\lambda}\Big]_{\lambda}=\Big\{\frac{1}{4}\cdot a_{\lambda},\frac{3}{4}\cdot a_{\lambda}\Big\}=\Big[\frac{1}{4}\cdot a_{\lambda}\Big]_{\lambda}\in L_{a_{\lambda}}.
$$
Note that $L_{a_{\lambda}}$ is an arc in $J_{\lambda}$ with end-points $[(0,0)]_{\lambda}$ and $[\frac{1}{2}\cdot a_{\lambda}]_{\lambda}$. Therefore, $q_{\lambda}\Big(\frac{3}{4}\cdot a_{\lambda}\Big)$ is a point in the interior of the arc $L_{a_{\lambda}}$. Next, for each positive integer $k$, let $e_k\in E(F_{\lambda})\setminus\{a_{\lambda}\}$ be such that 
    $$
    x_k\in K_{e_k}\setminus \{(0,0),e_k\},
    $$
 let $s_k\in (0,1)$ be  such that $x_k=s_k\cdot e_k$ and let 
$$
A_k=\{t\cdot e_k \ | \ t\in [0,s_k]\}  ~~~ \textup{ and } ~~~  B_k=\{[t\cdot e_k]_{\lambda} \ | \ t\in [0,s_k]\}.
$$
Also, let
$$
A=\Big\{t\cdot a_{\lambda} \ | \ t\in \Big[0,\frac{3}{4}\Big]\Big\}  ~~~ \textup{ and } ~~~  B=\Big\{[t\cdot a_{\lambda}]_{\lambda} \ | \ t\in \Big[0,\frac{1}{4}\Big]\Big\}.
$$
Note that for each positive integer $k$, $x_k$ is an end-point of $A_k$, and that $\frac{3}{4}\cdot a_{\lambda}$ is an end-point of $A$. Since $\displaystyle\lim_{k\to \infty}x_k=\frac{3}{4}\cdot a_{\lambda}$ and since $F_{\lambda}$ is smooth, it follows that
$$
\lim_{k\to \infty}A_k=A.
$$
On the other hand, note that
\begin{enumerate}
	\item for each positive integer $k$, $q_{\lambda}(x_k)$ is an end-point of $B_k$, 
	\item $q_{\lambda}(\frac{3}{4}\cdot a_{\lambda})=[\frac{1}{4}\cdot a_{\lambda}]_{\lambda}$ is an end-point of $B$, 
	\item $\displaystyle\lim_{k\to \infty}q_{\lambda}(x_k)=q_{\lambda}\Big(\frac{3}{4}\cdot a_{\lambda}\Big)$, and
	\item $\displaystyle\lim_{k\to \infty}B_k=L_{a_{\lambda}}  ~~~ \textup{ but } ~~~  L_{a_{\lambda}}\neq B$. 
\end{enumerate} 
This proves that the fan $J_{\lambda}$ is not smooth.
\end{proof}
\begin{theorem}\label{noth}
	For all $\lambda_1,\lambda_2\in \Lambda$, 
	$$
	\lambda_1\neq \lambda_2 ~~~ \Longrightarrow ~~~ F_{\lambda_1} \textup{ and } F_{\lambda_2} \textup{ are not homeomorphic.}
	$$
\end{theorem}
\begin{proof}
Let $\lambda_1,\lambda_2\in \Lambda$ be such that $\lambda_1\neq \lambda_2$. Then $F_{\lambda_1}$ is not homeomorphic to $F_{\lambda_2}$ since there is $e\in E(F_{\lambda_1})\setminus \{a_{\lambda_1}\}$ such that for each $f\in E(F_{\lambda_2})\}$,
	$$
	 |A_{F_{\lambda_2}}[(0,0),f]\cap\textup{JuMa}(F_{\lambda_2})|\neq |A_{F_{\lambda_1}}[(0,0),e]\cap \textup{JuMa}(F_{\lambda_1})|.
	$$
	Choose and fix such a point $e\in E(F_{\lambda_1})\setminus \{a_{\lambda_1}\}$. Note that  
	$$
	A_{F_{\lambda_1}}[(0,0),a_{\lambda_1}]\cap \textup{JuMa}(F_{\lambda_1})=\emptyset.
	$$
	 It follows from the definition of the relations $\sim_{\lambda_1}$ and $\sim_{\lambda_2}$ that for each element  $[f]_{\lambda_2}\in E(J_{\lambda_2})$,
	$$
	 |A_{J_{\lambda_2}}[[(0,0)]_{\lambda_2},[f]_{\lambda_2}]\cap\textup{JuMa}(J_{\lambda_2})|\neq |A_{J_{\lambda_1}}[[(0,0)]_{\lambda_{1}},[e]_{\lambda_1}]\cap \textup{JuMa}(J_{\lambda_1})|.
	$$
	Therefore, by Corollary \ref{jumas}, $J_{\lambda_1}$ and $J_{\lambda_2}$ are not homeomorphic. 
\end{proof}
Finally, we prove the main result of the paper.
\begin{theorem}
	There is a family $\mathcal H$ of uncountably many pairwise non-homeo\-morphic fans that are not smooth and that admit transitive homeomorphisms. 
	\end{theorem}
\begin{proof}
Let 
$$
\mathcal H=\{J_{\lambda} \ | \ \lambda \in \lambda\}.
$$ 
It follows from Theorem \ref{noth} that $|\mathcal H|=|\mathcal E|$. Since the family $\mathcal F$ from \cite[Definition 4.50]{banic1} is uncountable and since $|\mathcal E|=|\mathcal F|$, it follows from Theorem \ref{smoothf} that $\mathcal H$ is a family of uncountably many pairwise non-homeo\-morphic non-smooth fans. By Theorem \ref{trf}, for each $X\in \mathcal H$, $X$ admits a transitive homeomorphism.  
\end{proof}
\section{Acknowledgement}
This work is supported in part by the Slovenian Research Agency (research projects J1-4632, BI-HR/23-24-011, BI-US/22-24-086 and BI-US/22-24-094, and research program P1-0285). 
	

\noindent I. Bani\v c\\
              (1) Faculty of Natural Sciences and Mathematics, University of Maribor, Koro\v{s}ka 160, SI-2000 Maribor,
   Slovenia; \\(2) Institute of Mathematics, Physics and Mechanics, Jadranska 19, SI-1000 Ljubljana, 
   Slovenia; \\(3) Andrej Maru\v si\v c Institute, University of Primorska, Muzejski trg 2, SI-6000 Koper,
   Slovenia\\
             {iztok.banic@um.si}           
     
				\-
				

                 	\-
					
  \noindent J.  Kennedy\\
             Department of Mathematics,  Lamar University, 200 Lucas Building, P.O. Box 10047, Beaumont, Texas 77710 USA\\
{{kennedy9905@gmail.com}       }    

	\-
				
		\noindent C.  Mouron\\
             Rhodes College,  2000 North Parkway, Memphis, Tennessee 38112  USA\\\
{{mouronc@rhodes.edu}       }    

                 	\-
				
		\noindent V.  Nall\\
             Department of Mathematics,  University of Richmond, Richmond, Virginia 23173 USA\\
{{vnall@richmond.edu}       }   



\end{document}